%% file: SYMPLECTIC1.tex
\documentclass[11pt]{article}
 \pdfoutput=1
\usepackage{tikz-cd}
\usepackage{amsmath, amsthm, amsfonts}
\usepackage{xcolor}
\newtheorem{thm}{Theorem}[section]
\newtheorem{cor}[thm]{Corollary}
\newtheorem{lem}[thm]{Lemma}
\newtheorem{prop}[thm]{Proposition}
\theoremstyle{definition}
\newtheorem{defn}[thm]{Definition}
\theoremstyle{remark}
\newtheorem{rem}[thm]{Remark}

\def\cF{\mathcal{F}}

\date{\vspace{-5ex}}
\title{Symplectic and Poisson geometry of the information geometry’s Frobenius manifold.}
\author{Noemie Combe, Philippe Combe, Hanna Nencka \\
  \small Max Planck Institute for Mathematics in Sciences\\
  \small Inselstr. 22,\\
  \small 04103, Leipzig
 }

\begin{document}
\maketitle
{\bf keywords:} Frobenius manifolds; Symplectic geometry; Poisson structures;
\smallskip 

\begin{abstract}
We prove that the information geometry's Frobenius manifold is a symplectic manifold having Poisson structures. By proving this statement, a bridge is created between the theories developed by Vinberg, Souriau and Koszul and the Frobenius manifold approach. Until now these approaches seemed disconnected.  
 \end{abstract}

\section{Introduction} 
A fruit of the remarkable interaction between the physics vision of Quantum Field theory and the mathematical branches of algebraic geometry and algebraic topology was the creation of  Frobenius manifolds and Frobenius algebras. These mathematical objects entered naturally the scene in the process of axiomatisation of Topological Field Theory (TFT). 

\smallskip 

Until 2019, only three main classes of Frobenius manifolds had been explicitly known. Among these classes one has the $A$-side and $B$-side of the mathematical interpretation of the mirror problem.  In~\cite{CoMa20} an interesting connection between Frobenius manifolds and information geometry was shown. 
It turned out that there exists a fourth class of Frobenius manifolds, which is the class of manifolds of probability distributions, related to exponential families and having flat connection. Along this paper, this manifold will go by the name of the {\it ``fourth Frobenius manifold''} or the {\it ``information geometry's  Frobenius manifold''}. The ``fourth Frobenius manifold'', corresponds to the list in the classification of Frobenius manifolds outlined by Manin in~\cite{MaF}.

\smallskip 
The topic of this article is to prove that information geometry's Frobenius manifold is a {\it symplectic manifold} and that it forms a {\it Poisson manifold}. 
Consequences of such a result are two fold. First of all, this is a first step towards understanding a natural question, relating Frobenius manifolds and information geometry. Indeed, since Frobenius manifolds are a structure which appeared in the process of axiomatisation of Topological Field Theory TFT, the natural question is to understand the relation between the two dimensional (TFT) and the domain of information geometry. 

\smallskip 

On the second hand, we establish a direct bridge between the theory developed by Souriau--Koszul and the aspect given by the fourth Frobenius manifold. Indeed, initially Frobenius manifolds were defined using equations of hydrodynamical type. This approach was  developed by Dubrovin--Novikov (see for instance~\cite{DuNo}). Such an equation of hydrodynamical type implies the existence of Poisson structures, related to the fourth class of Frobenius manifolds. We put in parallel, different aspects: the Vinberg cone theory, Koszul--Souriau's approach along with ours. By doing so we establish a bridge in relation to information geometry, between the theory developed by Vinberg, Souriau and Koszul (see for instance~\cite{Ba16,Ba20,Ko,Sou65} and the Frobenius manifold~\cite{BaNo,DuNo,Du96,MaF,Man99} aspect which, until now 
seemed disconnected.  

A Frobenius manifold is a manifold equipped with an algebraic structure. This algebraic structure is resumed in the two following axioms, known as the {\it associativity} and {\it potentiality} axioms. 
 \begin{itemize}
\item[$\ast$] {\bf Associativity:} For any triplet of (flat) local tangent fields denoted by $u,v,w$, the following is satisfied: \begin{equation}\label{Asso}t(u,v,w)=g(u\circ v,w)= g(u, v\circ w),\end{equation} where: 
\begin{itemize}
\item $\circ$ is a multiplication on the tangent sheaf.
\item $t$ corresponds to a tensor of rank three, 
\item $g$ is a bilinear form, being the metric.
\end{itemize}
\item[$\ast$] {\bf Potentiality:} $t$  admits locally everywhere a potential function $\varPhi$. For any local tangent fields $\partial_i$, we have the following: 
\begin{equation}\label{Pot}t(\partial_a,\partial_b,\partial_c)=\partial_a\partial_b\partial_c\varPhi.\end{equation}
\end{itemize}
The associativity condition can also be encapsulated in the Witten--Dijkgraaf--Verlinde--Verlinde (WDVV) non-linear PDE system (see \cite{Man99}, p. 19-20):
\begin{equation}\label{E:WDVV}\forall a,b,c,d: \underset{ef}{\sum}\varPhi_{abe}g^{ef}\varPhi_{fcd}=(-1)^{a(b+c)}\underset{ef}{\sum}\varPhi_{bce}g^{ef}\varPhi_{fad}.\end{equation}

The paper is devoted to proving the following main theorem:

\begin{thm}[Main theorem] 
The information geometry's Frobenius manifold is a symplectic manifold and has Poisson structures.
\end{thm}
From our perspective, we interpret the information geometry $F$-manifold as a certain type of configuration space with marked points in the Euclidean space.

The plan of the paper goes as follows:
\begin{itemize}
\item First of all, we prove the main theorem. This is split in two subsections within section 2, for the sake of clarity. One part aims at proving that we have an integrable system and the other one aims at showing the Hamiltonian structure.
\item We discuss the Hamiltonian action in section 2.3. The combination of both statements leads to having a symplectic manifold.
\item Section 3 deals with connections between theories developed by Vinberg, Souriau and Koszul. Namely, this is about the Poisson structures and their interaction with the Frobenius manifold structure. There is a deep relation between Poisson structures of hydrodynamical type and Frobenius algebras. 
\end{itemize}

\input{S1IntegrableSystem}

\input{PoissonBrackets}

\bibliographystyle{alpha}

\bibliography{Biblio}

\end{document}

%% file: S1IntegrableSystem.tex
\section{Proof of the main theorem}

We expose the proof of the main theorem. It is known that a symplectic manifold is an integrable system, equipped with a Hamiltonian action (for sinatance see \cite{Mo09}), relying on the metric tensor of the manifold. As was presented in the introduction, the discussion is decomposed into two principal steps:

\begin{itemize}
\item The fourth Frobenius manifold is an integrable system.
\item This integrable system is a Hamiltonian one, depending on the tensor metric of the manifold. Smaller steps are needed to prove this statement, (we choose not mention them here, for simplicity). Recollections on Hamiltonian systems and paracomplex geometry are needed for the reader's convenience. 
\end{itemize}
\medskip 
\subsection{First step: integrable system}
 To prove that we have an integrable system, we invoke the following theorem of Mokhov~\cite{Mo07,Mo09} stating that:

\smallskip 

{\it The class of flat, torsionless submanifolds in any Euclidean or pseudo-Euclidean space is an integrable system.}

\smallskip

Now, we know from \cite{Du96,CoMa20,Man99} that the Frobenius manifold is necessarily {\it flat and torsionless}. Therefore, the information geometry fourth Frobenius manifold is an integrable system. This proves the first step. 

\medskip 

\subsection{Second step: Hamiltonian system} The second steps means that we ought to show that the system is a Hamiltonian one.
Namely, the important statement to use is that  the system {\it has} a Hamiltonian structure if and only if the {\it Lagrangian $\mathcal{L}$ depends on the tensor metric}, i.e. $\mathcal{L}(x,\xi)=g_{ij}(x)\xi^i(x)\xi^j(x) - U(x),$
where $x$ are canonical coordinates on the manifold, $g_{ij}$ is the tensor metric, $U$ is a scalar field and $\xi^i,\xi^j$ are vector fields. 
\smallskip


\input{Hamiltonian}
\smallskip 
\subsubsection*{Hamiltonian construction for the fourth Frobenius manifold}
We continue by showing the explicit construction of the Hamiltonian. From~\cite{CoMa20}, it is known that the information geometry $F$-manifold has a {\it paracomplex structure}. It implies that the construction of differential forms for paracomplex  geometry is necessary. Therefore, we recall a few notions on paracomplex geometry. 

\smallskip 
\subsubsection*{Recollections on paracomplex geometry}
Let us define the rank 2 algebra of paracomplex numbers $\frak{C} = \mathbb{R} + e\mathbb{R}$. This is a split algebra, generated by the basis $\{1,e\}$ where $e^2=1$. It resembles the field of complex numbers in the sense that a paracomplex number $z=x+ey$ has a paracomplex conjugate $\overline{z}=x-ey$ where $e$ squares to +1. However, this algebraic structure does not form a field. On the other side, an interesting feature of this split algebra is that it has a {\it pair of idempotents}.

An idempotent $a$ in an algebra is element such that $a^2=a.$ In the algebra of the paracomplex numbers the pair of idempotents can be seen using a suitable change of basis. The new basis is generated by $e_{+}$ and $e_{-}$ where $e_{\pm}=\frac{1\pm e}{2}$. It can be easily checked that $e_+^2=e_+$ (resp. $e_-^2=e_-$). This pair of idempotents plays an important geometric role, since it implies the existence of a Peirce reflection (also known as Peirce mirror).   

\smallskip 

Let us roughly recall the notions of paraholomorphic structures and manifolds.  Denote by $E_{2m}$ a $2m$-dimensional real affine space. {\it A paracomplex structure} on $E_{2m}$ is an endomorphism $\frak{K}: E_{2m} \to E_{2m}$ such that $\frak{K}^2=I$,
and the eigenspaces $E_{2m}^+, E_{2m}^-$ of $\frak{K}$ with eigenvalues $1,-1$ respectively, have the same dimension. A {\it paracomplex manifold} is a real manifold $M$, equipped with a paracomplex structure $\frak{K}$, admitting an atlas of paraholomorphic coordinates---being functions with values in the split algebra of paracomplex numbers---and such that the transition functions are paraholomorphic.

\smallskip 

Considering a module over the paracomplex algebra, and more precisely its (real) realisation as a linear space, a domain $D$ in the paracomplex space contains paracomplex points $z=(z^1,\cdots,z^n)$, which are given by:

\[z^k=x^k+e y^k,\]

and by definition we have $e^2=1$ and $k\in \{1,\cdots,n\}$. 

\smallskip 

The existence of the paracomplex structure implies the existence of local coordinates $(z_+^\alpha, z_-^\alpha),\, \alpha = 1,\cdots, m$ such that
paracomplex decomposition of the local tangent fields is of the form
\[
T^{+}M=span \left\{ \frac{\partial}{\partial z_{+}^{\alpha}},\, \alpha =1,\cdots, m\right\},
\]
\[
T^{-}M=span \left\{\frac{\partial}{\partial z_{-}^{\alpha}}\, ,\, \alpha =1,\cdots, m\right\} .
\]
Such coordinates go by the name of {\it adapted coordinates} for the paracomplex structure $\frak{K}$.

The dual basis for the differential forms (dual basis to the paracomplex tangent vectors) is given by:

\[dz^k=dx^k+edy^k, \quad d\overline{z}^k=dx^k-edy^k,\quad k=\{1,\cdots,n\}.\] 

\smallskip 

Similarly to the complex case, we have a splitting of the differential forms:
\[d=d'+d'',\] where for a given $(p,q)$ paracomplex form $\omega$ the differential $d'\omega$ is a $(p+1,q)$-form, whereas $d''\omega$ is a$(p,q+1)$-form. It is easily seen that \[(d')^2=(d'')^2=0\] and that \[d''d'=-d'd''.\]

\smallskip 

In order to obtain the Hamiltonian on the paracomplex manifold, we introduce a local paracomplex Dolbeault (1,1)-form. This is given by:

\[\omega:=\partial_+\partial_-\phi\]
where $\phi$ is the potential function defined uniquely modulo subspace of local functions $ker \partial_+\partial_-$ and $\partial_\pm$ are given in canonical paracomplex coordinates. 

To show that the Hamiltonian relies on the metric tensor, it is sufficient to consider the scalar product on a space with paracomplex structure. By definition the scalar product is defined as:
\[\langle \xi, \eta \rangle_{\frak{C}}= g_{jk}\xi^j\overline{\eta^{k}}.\] 
This scalar product is a Hermitian one, since we have the equality:
\[\langle \xi, \eta \rangle_{\frak{C}}=\overline{\langle \eta,\xi \rangle_{\frak{C}}}.\]

Then, the paracomplex differentiable 2-form is given by: \[\Omega_{\frak{C}}=\frac{e}{2}g_{jk}dz^j\wedge d\overline{z}^k.\]

By definition, it is known that the Hamiltonian is a 2-form, so we have proved that a Frobenius manifold---with paracomplex structure---is an integrable system with Hamiltonian.  

\medskip 
\subsection{The symplectic fourth Frobenius manifold}
It remains to adapt this construction to the information geometry Frobenius manifold. 


In order to show that this manifold is symplectic it is sufficient to notice that the Hamiltonian discussed previously is associated to the following functional:

\begin{equation}\label{E:1}
\langle \mu, f\rangle=\int_{\tilde{\Omega}}f(\omega)\mu\{d\omega\},\end{equation}
where $f$ is a $\cF$-measurable function and $\mu$ is a dominated measure on the measurable space $(\tilde{\Omega},\cF)$. In the discrete case, this is given by: 

\[\langle \mu, f\rangle=\sum_{j=1}^m f^j\mu_j.\]



We want to express the Hamiltonian for statistical manifolds, being Frobenius. From \cite{CoMa20,CoCoNen21} we have that the statistical manifolds (satisfying the Frobeniusity condition) have geometry of the Lie group of type $SO(1,n-1)$. The bilinear form for this group, is given by :
\[\kappa_1ds^2=\lambda_{\mu\nu} dz^{\mu}dz^{\nu},\]
where $\{\lambda_{\mu\nu} \}$ is a diagonal matrix $(1,1,1,-1)$ and $\kappa_1$ is a constant.

So, the Lagrangian written in terms of the tangent vector is given by: $\mathcal{L}=\frac{1}{2}C(\xi^\mu\xi_{\mu}-1)+ \kappa_2\xi^\mu A_{\mu}$
where $A_{\mu}$ corresponds to the Gauge Abelian field; $C$ and $\kappa_2$ are constants. 
The so-called ``momentum'' and ``force'' are derived from the Lagrangian, respectively as follows:
\[p_{\mu}=\frac{\partial \mathcal{L}}{\partial \xi^\mu}=\xi_{\mu}+\kappa_2A_{\mu},\]
\[f_{\mu}=\frac{\partial \mathcal{L}}{\partial z^\mu}=\xi^\mu\frac{\partial A_{\nu}}{\partial z^\mu}.\]

This allows to write the Hamiltonian using the Lagrangian view point for our case. In other words we have:
\[\mathcal{H}=\xi^\mu\frac{\partial\mathcal{L}}{\partial \xi^\mu}-\mathcal{L}=\xi^\mu p_{\mu}-\mathcal{L}.\]
So, in this way we show that our information geometry manifold has Hamiltonian structure.    
Finally, since the manifold under consideration is Frobenius, this implies it has torsionless and flatness properties. Therefore, by the above statement it forms an integrable system. Now, the existence of the Hamiltonian implies that this Frobenius manifold is symplectic. $\square$

\smallskip 

We have thus ended the proof of the main theorem. In the next part, we discuss other implications of this on the geometry.

%% file: Hamiltonian.tex
\subsubsection{Recollections on the construction of the Hamiltonian system}
For the convenience of the reader, some recollections are presented (in this subsection) in relation to the Hamiltonian system.
Let $M$ be a manifold of dimension $n$ and let $x=(x^1,...,x^n)$ be the canonical coordinates attached to a point $x$ on $M$. 
We construct a Hamiltonian system using the $2n$-dimensional fibered cotangent bundle $T^*(M)$.
Given a point $x\in M$ and a covector $p=(p_1,...,p_n)$ attached to $x$, we have the following closed (and exact) differential form given by:
\[ \sum_{\alpha}dx^{\alpha}\wedge dp_{\alpha}.\]

\smallskip 

Note that the differential form $\Omega=d\omega,$ where $\omega$ is given by $p_{\alpha}dx^\alpha$ is defined globally on $T^*(M)$ if and only if the form is {\it proportional} to the volume element (with a non-nul proportionality coefficient). In particular, this form defines the scalar product (alternate and non-degenerate) of a pair of vectors:
\[\langle \xi, \eta \rangle=-\langle \eta, \xi\rangle=J_{ij}\xi^i\eta^j \]
where $i,j=1,2,...,2n;$\, $\Omega=J_{ij}dy^i\wedge dy^j$, with non-nul $J_{ij}$ and $(y^1,\dots,y^{2n})$ are canonical coordinates on $T^*(M)$ given by:

\[y^i=x^i \quad y^{n+i}=p_{i},\quad i=1,2,...,n.\]
Using this, we establish a relation to phase space (in the physics setting) and therefore to configuration spaces of marked points in the euclidean space.

Now, this leads to the Hamilton equation:

\[\frac{dy}{ds}=J^{ij}\frac{\partial \mathcal{H}}{\partial y^j}, \quad J^{ij}J_{jk}=\delta_{k}^{i}.\]
where $s$ is a real parameter.


We present the explicit construction of the Hamiltonian for the case of the information geometry's $F$-manifold. 
As it is known, the energy $\mathcal{E}$ is given by $\mathcal{E} = \mathcal{H}(x, p)$ where  
\[\mathcal{H}(x, p) = \frac{1}{2} \langle p\, ,\, p \rangle,\]
and $p$ was defined above. 

Since $p$ can be considered as a linear operator, given by $\frac {\partial }{\partial x}$
then the Hamiltonian can be written as a bi-linear form:
 $\mathcal{H} (x\, ,\, \frac{\partial}{\partial x})=\mathcal{H}(x\, ,\, \xi)$. This corresponds to the 2- form: $\langle\frac {\partial}{\partial x}\, ,\, \frac{\partial}{\partial x}\rangle=\langle\xi\, ,\, \xi\rangle.$
In the following we show an explicit Hamiltonian for the information geometry fourth Frobenius manifold. 

%% file: PoissonBrackets.tex
\section{Poisson structures}
In this section, we present the Poisson bracket for the information geometry $F$-manifold. We consider a symplectic manifold. Novikov, in his works on symplectic geometry, introduced a special type of Poisson brackets for this type of manifold. More specifically, we called them the {\it Novikov--Poisson brackets}. Initially, these Novikov--Poisson brackets were introduced to construct a theory of conservative systems of hydrodynamic type. It strongly interferes in the existence of Frobenius manifold structure.

By the previous section, we have the following  proposition:

 \begin{prop}
The information geometry fourth Frobenius manifold  is a symplectic (paracomplex) manifold $(M,\Omega_{\frak{C}},\Omega)$ with the 2-form $\Omega_{\frak{C}}=\frac{e}{2}g_{jk}dz^j\wedge d\overline{z}^k,$ where $e^2=1$ and $g_{jk}$ is a paracomplex bilinear form and $\Omega=dz^j\wedge dp_k.$
\end{prop}

This can be enlarged to Poisson structures and Poisson manifolds.  A Poisson structure is a Lie algebra bracket $\{\cdot\,,\, \cdot\}$ on the vector space of smooth functions on a manifold $M$ which form a Poisson algebra with the point-wise multiplication of functions. Therefore, we can say that:
\begin{cor}
The information geometry fourth Frobenius manifold carries a Poisson structure and thus is a Poisson manifold. 
\end{cor}

The Poisson bracket on the paracomplex space $\{\cdot\, ,\,\cdot \}_{\frak{C}}$ is given by \[\{\xi,\eta\}_\frak{C}=\frac{1}{2}Im \langle \xi,\eta \rangle_{\frak{C}}=\frac{e}{2}\left(\frac{\langle \xi, \eta \rangle_{\frak{C}}- \overline{\langle \xi, \eta \rangle_{\frak{C}}}}{2}\right),\] where $e$ satisfies $e^2=1$ and $Im$ is the imaginary part of a paracomplex number. 

We consider this Frobenius manifold as a certain type of phase space, i.e. $(z^1,\cdots, z^n,\, p_1,\cdots, p_n, \Lambda_1,\cdots, \Lambda_m)$,
where the $z^\mu$ are paracomplex numbers, $p_\mu=p_\mu(z)$ are 1-forms of the fibered cotangent bundle attached at the point $z$ and the $\Lambda_i$ are $m$ functions on the phase space. For instance, concerning the set of functions $\{\Lambda_i\}_{i=1}^m$ on the phase space we can consider them as a set of spins of Ising, Heisenberg or any other relevant structure. In the case of the Ising spin model we have $m=1$. Whereas, for the case of the Heisenberg model $m=3$.

\smallskip 

The set $\{\Lambda_i\}_{i=1}^m$ can have a Lie algebra structure, if the Poisson bracket of a pair of functions $\Lambda_i$ and $\Lambda_j$ is obtained by a linear combination of those functions, for $i,j\in \{1,\cdots,m\}$. More precisely, we have $\{\Lambda_i,\Lambda_j\}=-\Lambda_k\gamma^k_{ij}$, where $\gamma^k_{ij}$  are structure constants of the algebra. 

Let $A, B$ be a pair of functions on the phase space $(z^1,\cdots, z^n,p_1\cdots, p_n)$. These are observables and can be for example the momentum. $A$ and $B$ are real functions depending on $z$ and $p$. Then, the Poisson bracket is given by:
\[\{A,B\}=\frac{\partial A}{\partial p_\mu}\frac{\partial B}{\partial z^\mu}-\frac{\partial B}{\partial p_\mu}\frac{\partial A}{\partial z^\mu}.\]

The Poisson bracket can be interpreted as a type of derivation on the functions $A$ and $B$ in the sense that for any given pair of smooth phase functionals $f(A)$ and $g(B)$ we have the following relation:

\[\{f(A),g(B)\}=f'(A)g'(B)\{A,B\}.\]
In the generalised context, i.e. for $A=(A_1,\cdots, A_r)$ and $B=(B_1,\cdots B_l)$ we have:
\[\{f(A),g(B)\}=\frac{\partial f}{\partial A_i}\frac{\partial g}{\partial B_j}\{A_i,B_j\}.\]

\begin{rem}
Note that the derivation of smooth functions correspond to smooth tangent vector fields (for every given smooth function). So, for a smooth function $f$ there exists a vector $X_f$ given by $X_f(g)=\{f,g\}$. This is also known as the {\it Hamiltonian vector field}, corresponding to the function $f$. 
\end{rem}

In the most general setting, that means in the generalised phase space where the $\Lambda_i$ functions enter the game, $(z^1,\cdots, z^n,\, p_1,\cdots, p_n, \Lambda_1,\cdots, \Lambda_m)$ we have the Poisson bracket defined as:
\[\{A,B\}=\frac{\partial A}{\partial p_\mu}\frac{\partial B}{\partial z^\mu}-\frac{\partial B}{\partial p_\mu}\frac{\partial A}{\partial z^\mu}-\Lambda_k\gamma^k_{ij}\frac{\partial A}{\partial \Lambda_i}\frac{\partial B}{\partial\Lambda_j}.\] 

Note that for an arbitrary Hamiltonian $\mathcal{H}(p,z,\Lambda,s)$ and a smooth phase space function $Q(p,z,\Lambda)$ the equation of motion is described by $\dot{Q}=\frac{dQ}{ds}=\{\mathcal{H},Q\}.$ This is identified to the notion of a {\it symplectic gradient}. Indeed, a symplectic gradient is
a unique vector field $X_{\mathcal{H}}$ described by the condition $X_{\mathcal{H}}(Q)=\{\mathcal{H},Q\}$. For instance, one can write the following symplectic gradient: $X_{\mathcal{H}}(z^\mu)=\frac{d z^\mu}{ds}= \{\mathcal{H},z^\mu\}$ where $\mathcal{H}$ is the Hamiltonian, mentioned and constructed in the above section. More formally, the symplectic gradient of $\mathcal{H}$ is the vector field given by:
\[X_{\mathcal{H}}:=\Omega_{\frak{C}}^{-1}d_{dR}\mathcal{H},\]
where the symbol $d_{dR}$ is the {\it de Rham differentia}l, mapping smooth functions on $M$ to 1-forms on $M$. On the fiber bundle space $T^*(M)$, we can indeed define a symplectic gradient. Since the 2-form $\Omega_{\frak{C}}$ is closed, the Poisson brackets introduce on the linear space of functions a Lie algebra structure. 



We have thus the important corollary,  that the fourth Frobenius manifold is a Poisson manifold and is thus extensively equipped with Poisson Geometry. 

\smallskip 

The following discussion, tends to show a connection of our latter statement with the works in~\cite{Ba20}. In particular, our last statement establishes a bridge to the works of Koszul--Souriau \cite{Ba16,Ba20,Ko,Sou65}.

The Frobenius structure is closely related to the Poisson bracket of hydrodynamic type (see \cite{BaNo,DuNo}). More precisely, we present the definition of Poisson bracket of hydrodynamic type and discuss the consequences. 

We define the Novikov-type of Poisson bracket, based on the works of~\cite{BaNo,DuNo}. 

\begin{defn} The Novikov--Poisson bracket is a Poisson bracket of hydrodynamic type, defined by the following equation:
 \begin{equation}\label{E:p}
  \{u^i(x),u^j(y)\} = g^{ki,\alpha}(u(x))\partial_\alpha\delta'(x - y) + b^{ki,\alpha}_j(u(x))\partial_\alpha u^j\delta(x - y),\end{equation} 
  where $u^j$ are fields (usually the density of momentum and energy or mass) and $b^{ij}_k$ is a tensor. 
\smallskip 
 \end{defn}
 
Remark the following. 
\begin{lem}
Let $z=(z^1,\cdots,z^n)$ and $w=(w^1,\cdots,w^n)$ be two points on a paracomplex manifold and let $u(z)=(u^1(z),\cdots, u^n(z))$, $v(w)=(v^1(w),\cdots, v^n(w))$ (with $u^i(z)=u(z^i)$) be two ''observables'' on the phase space, associated to the paracomplex Frobenius manifold. 
Then, $\{u^i(z)\, ,\, u^j(w) \}=g^{ij}(u(z)\delta'(z-w))+\frac{\partial u^k}{\partial z}b^{ij}_{k}(u(z)\delta(z-w))$.
\end{lem}
\begin{proof}This follows from the discussion in section 2 on paracomplex spaces, from the result in\cite{CoMa20} showing the existence of paracomplex geometry of the fourth Frobenius manifold and from \cite{DuNo}.\end{proof}
\smallskip 

The classification of Poisson brackets has been done in (\cite{BaNo}). This classification depends linearly on $u^k$ and is relative to linear changes $u^k = A^k_j w^j$. It can be obtained through the classification of Lie algebras and Frobenius algebras. Indeed, the {\it simplest local Lie algebra} arises from the brackets as follows: 

\[g^{ij}=C^{ij}_ku^k+g_{0}^{ij}, \quad b_k^{ij}=\text{cst}, \quad g_{0}^{ij}=\text{cst};\]
\[ [p,q]_k(z) = b^{ij}_k(p_i(z)q_j'(z)- q_i(z)p'_j(z))\]
 where 
 \begin{equation}\label{E:b}
  b^{ij}_k+b^{ji}_k=C^{ij}_k=\partial g^{ij}/\partial u^k.
  \end{equation} 
 \begin{rem}
It is precisely the tensor $b^{ij}_k$ (from equation\, \eqref{E:b}) which plays an important role from the side of the Lie algebra. Indeed, it defines a local, translationally invariant Lie algebra of first order, if and only if the following condition is fulfilled. 
\begin{center}
{\bf Condition:} {\it The tensor $b^{ij}_k$ plays the role of the constant structure in the multiplication law of a finite dimensional algebra $\mathcal{A}$ defined by the identities:}
 \smallskip 
 
\[a,b,c\in \mathcal{A},\quad e^ie^j\quad e^ie^j=b_k^{ij}e^k;\]
\[a(bc)=b(ac),\quad (ab)c-a(bc)=(ac)b-a(cb).\] 

\end{center}

The algebra $\mathcal{A}$ is associative and commutative if the condition $2b_{ij} = 2b_{ji} = C_{ij}$ is fulfilled. By definition, this algebra is called Frobenius if there exists a given  non-degenerate inner product satisfying the associativity condition i.e. such that: $\langle e_i\circ e_j\, ,\, e_k \rangle=\langle e_i\, ,\, e_j\circ e_k \rangle.$ The necessary
and sufficient condition for an algebra with identity such as $\mathcal{A}$ to be Frobenius is that $2b_{ij}u^k = C_{ij}u^k$ is {\it non-degenerate} at a generic point. 

\smallskip  
Let us explicit our thought. The $b^{ki,\alpha}_j$ (resp. $b_{ji}^k$ in a less general context) plays a decisive role, regarding the Frobenius structure in the following way. Those tensors correspond to the constant structures in the multiplication law of the Frobenius algebra $\mathcal{A}_p$, associated to the Frobenius manifold $M$. More precisely, the Frobenius algebra is associated to the tangent space of a Frobenius manifold $T_pM$, at a given point $p$. This leads to the formulation of the statement below.  

\begin{lem}
For the information geometry's Frobenius manifold, there exist Poisson brackets of hydrodynamic type, (equation \eqref{E:p}) such that the tensors $b^{ij}_k$ correspond to the constant structures in the multiplication law of the Frobenius algebra, associated to the tangent space to the Frobenius manifold.
\end{lem}
\begin{proof}
This follows from the definition of Frobenius manifold introduced by Dubrovin in~\cite{Du96} and from the principal statement in~\cite{BaNo}. 
\end{proof}
In particular, we have the following statement:

\begin{lem}
Let $(M,\Omega)$ be the symplectic information geometry Frobenius manifold. Then, the phase space associated to it has an algebraic structure depending on the tensors $b^{ki,\alpha}_j$, defined in equation \eqref{E:p}. 
\end{lem}
\begin{proof}
Let us map $T_pM$ to its dual space $T^*_pM$ by a linear operator $\Phi$. Then, there exists a map $\tilde{\phi}:\mathcal{A}_p\to \mathcal{A^*}_p$, where    $\mathcal{A^*}_p$ is the algebra attached to $T^*_pM$. Following the commutative diagram below, we have that $\tilde{\phi}=\Psi^{-1}\circ\Phi\circ\chi$ is a homomorphism.
\[
\begin{tikzcd}
\mathcal{A}_p \arrow{r}{\tilde{\phi}} \arrow[swap]{d}{\chi} & \mathcal{A^*}_p \arrow{d}{\Psi} \\
T_pM  \arrow{r}{\Phi} & T^*_pM
\end{tikzcd}
\]

Therefore, by duality (and only for finite dimensional cases!) the $b^{ki,\alpha}_j$ are the constant structures of the algebra associated to the cotangent space $\mathcal{A^*}_p$. So, the coefficients $b^{ki,\alpha}_j$ determine the underlying algebraic structure of the phase space and thus of the symplectic paracomplex manifold. 

\end{proof}

\smallskip 

Finally, one last remark concerning the Poisson geometry of this Frobenius manifolds. The  Novikov--Poisson brackets  (hydrodynamic type) are tightly related to the metric of the manifold in the following way. 

\begin{prop}
On the fourth Frobenius manifold, there exist flat pencils of Koszul--Souriau metrics, being related to compatible Novikov--Poisson brackets.
\end{prop}
\begin{proof}
Suppose we have a flat metric in some coordinate system $\{x^i\}$. If the components of the metric $g^{ij}(x)$ and the Christoffel symbols $\Gamma^{ij}_k$ of the corresponding Levi--Civita connection depends linearly on the coordinate $x^1$, then the metrics form a {\it flat pencil:} 

\[g^{ij}\quad \text{and}\quad  g_2^{ij}:=\partial g^{ij}\]

(it is assumed that the determinant of $g_2^{ij}$ is not zero). 

A pair of Poisson brackets  (hydrodynamic type) $\{\cdot\, ,\, \cdot\}$ are called {\it compatible} if an arbitrary linear combination with constant coefficients of those brackets defines back a Poisson bracket. Flat pencils of metrics correspond to {\it compatible pairs of poisson brackets} (hydrodynamic type) (see \cite{Du96}). Therefore, we have the existence of flat pencils of Koszul--Souriau metrics for the fourth Frobenius manifolds.
\end{proof}

\smallskip

\section{Remarks}
We can remark the following facts. 
\begin{rem} Orbits (that means geodesics) in the fourth Frobenius manifold are given by some trajectories given in paraholomorphic coordinates. 
Motions along these paraholomorphic orbits are {\it non dissipative}~\cite{NeSe}. In particular, by non dissipative motion we mean that  the production of entropy (during the motion) is null. In the latter case---and referring to Souriau's theory---we have an equality between the inner statistical states and the outer statistical states. Whereas the motion of transversal vectors to these surfaces is dissipative \cite{Ba20,BuNePe}.
\end{rem}

 \end{rem}
 \begin{rem} 
The orbit (that means a given geodesic) on a paracomplex space is decomposed into two orbits. Following the construction of paracomplex spaces defined in section 2, each of these orbits lies on a flat submanifold, corresponding respectively to one of the generators: $e_+$ or $e_-$ of the algebra. These two orbits are equidistanced to the {\it Peirce mirror}. 

\smallskip 

We digress on the notion of Peirce reflection (or Peirce mirror) for the convenience of the reader.  Peirce reflections arise for a given unital Jordan algebra (say $J$) with orthogonal idempotents $a, a'$ i.e. verifying the condition $a a'=a'a=0$. We are precisely in this case (see~\cite{CoCoNenFICC}).  
Considering a pair of orthogonal idempotents $a, a'$ and by putting $u=a-a'$ we can show that $u$ is an involution. Now, 
if $u$ is an involution in a unital Jordan algebra $J$, then one can define an involutionary automorphism of $J$, i.e. $U_u^2 = 1_J$ (see lemma 6.10 in \cite{MC04}). Then this involution has a description as a Peirce reflection in terms of the Peirce projections (see p.238 \cite{MC04} for more details). For further details and relations between the fourth Frobenius manifold and Peirce's reflection and decomposition theorems, see \cite{CoCoNenFICC}.

\smallskip 

In the light of the Peirce mirror, we can say that these orbits form a so-called pair of Clifford parallels. By Clifford parallel we mean a straight line in an elliptic space staying at constant distance from a given base (straight) line. The motion on paraholomorphic orbits is equivalent to the motion on pairs of Clifford parallels. Using the previous remark, we can state that the paraholomorphic motion of Clifford pairs is non dissipative. 
\end{rem}

